\documentclass[11pt]{amsart}
\usepackage{amssymb}
\usepackage{amsthm}
\usepackage{mathrsfs}
\usepackage{graphicx}
\usepackage{epstopdf}
\usepackage[all]{xy}
\newtheorem*{thmintro}{\bf Theorem}  

\newtheorem{thm}{Theorem}[section]
\newtheorem{cor}[thm]{Corollary}
\newtheorem{lem}[thm]{Lemma}
\newtheorem{prop}[thm]{Proposition}
\theoremstyle{definition}

\theoremstyle{remark}
\newtheorem{rem}[thm]{Remark}
\numberwithin{equation}{section}

\newcommand{\abs}[1]{\left\vert#1\right\vert}
\newcommand{\set}[1]{\left\{#1\right\}}

\newcommand{\Real}{\mathbb R}

\newcommand{\eps}{\varepsilon}

\newcommand{\sph}{\mathbb{S}}
\newcommand{\To}{\rightarrow}

\newcommand{\df}{\mathrm{d}}
\newcommand{\vphi}{\varphi}

\newcommand{\Hess}{\mathrm{Hess}}

\begin{document}

\title[Gradient Yamabe Solitons]{Gradient Yamabe Solitons on Warped Products}%
\author{Chenxu He}%
\address{Department of Mathematics, Lehigh University}%
\email{he.chenxu@lehigh.edu}%


\begin{abstract}
The special nature of gradient Yamabe soliton equation which was first observed by Cao-Sun-Zhang\cite{CaoSunZhang} shows that a complete gradient Yamabe soliton with non-constant potential function is either defined on the Euclidean space with rotational symmetry, or on the warped product of the real line with a manifold of constant scalar curvature. In this paper we consider the classification in the latter case. We show that a complete gradient steady Yamabe soliton on warped product is necessarily isometric to the Riemannian product. In the shrinking case, we show that there is a continuous family of complete gradient Yamabe shrinkers on warped products which are not isometric to the Riemannian product in dimension three and higher.
\end{abstract}
\maketitle

\section{Introduction}

Geometric flows are important tools to understand the topological and geometric structures in Riemannian geometry. A special class of solutions on which the metric evolves by dilations and diffeomorphisms plays an important role in the study of the singularities of the flows as they appear as possible singularity models. They are often called soliton solutions. In the case when the diffeomorphisms are generated by a gradient vector field, we call such soliton a gradient soliton. In this paper we are interested in the gradient soliton solutions to the Yamabe flow. This flow has been very well-understood in the compact case, see the very recent survey \cite{Brendle} by S. Brendle and the references therein. For the non-compact case, in \cite{DSesumYamabeflow} P. Daskalopoulos and N. Sesum showed that the solutions to the Yamabe flow from some complete metric develop a finite time singularity and the metric converges to a soliton solution after re-scaling.

A complete Riemannian manifold $(M, g)$ is called a \emph{gradient Yamabe soliton} if there exists a smooth function $f : M \rightarrow \Real$ and a constant $\rho \in \Real$ such that the Hessian of $f$ satisfies the equation
\begin{equation}\label{eqnGYamabeS}
\Hess f = \left(R - \rho\right)g,
\end{equation}
where $R$ is the scalar curvature of $M$. If $\rho > 0$, $\rho < 0$ or $\rho = 0$, then $(M, g)$ is called a \emph{Yamabe shrinker}, \emph{Yamabe expander} or \emph{Yamabe steady soliton} respectively. In \cite{DSesum} Daskalopoulos and Sesum initiated the investigation of gradient Yamabe solitons and showed that all complete locally conformally flat gradient Yamabe solitons with positive sectional curvature are rotationally symmetric. This result is inspired by the work of the classification of locally conformally flat Ricci solitons, especially in \cite{CaoChen}.

The equation $(\ref{eqnGYamabeS})$ where the right hand side is a smooth function, not necessarily the one given by the scalar curvature, was appeared in 1925 on the study of Einstein metrics that are conformal to each other by H. Brinkmann, see \cite{Brinkmann}. Such equation was also considered by J. Cheeger and T. Colding in their work \cite{CheegerColding}. By exploring the special nature of the Yamabe soliton equation (\ref{eqnGYamabeS}) H.-D. Cao, X. Sun and Y. Zhang showed that every complete gradient Yamabe soliton admits a global warped product structure in \cite{CaoSunZhang}.
\begin{thmintro}[Cao-Sun-Zhang]
Let $(M^n, g, f)$ be a complete gradient Yamabe soliton satisfying equation (\ref{eqnGYamabeS}) with non constant function $f$. Then $\abs{\nabla f}$ is constant on regular level hypersurfaces of $f$, and either
\begin{enumerate}
\item $f$ has a unique critical point, and $(M, g)$ is rotationally symmetric and equal to the warped product
\[
\left([0, \infty), \df r^2\right)\times_{\abs{\nabla f}}\left(\sph^{n-1}, \bar{g}_{\text{can}}\right),
\]
where $\bar{g}_{\text{can}}$ is the round metric on $\sph^{n-1}$, or
\item $f$ has no critical point and $(M, g)$ is the warped product
\[
\left(\Real, \df r^2\right) \times_{\abs{\nabla f}} (N^{n-1}, \bar{g}),
\]
where $(N^{n-1}, g)$ is a Riemannian manifold of constant scalar curvature, say $\bar{R}$. Moreover, if the scalar curvature $R$ of $(M, g)$ is non-negative, then either $\bar{R} > 0$, or $R = \bar{R} = 0$ and $(M, g)$ isometric to the Riemannian product $(\Real, \df r^2)\times (N^{n-1}, \bar{g})$.

\end{enumerate}
\end{thmintro}

By the global warped product structure from the theorem above, the Yamabe soliton equation (\ref{eqnGYamabeS}) reduces to the following ordinary differential equation in $\vphi = \abs{\nabla f}$, see \cite[Equation (2.17)]{CaoSunZhang}:
\begin{equation}\label{eqnphi2nd}
\vphi' + \rho = \vphi^{-2}\bar{R} - (n-1)(n-2)\left(\frac{\vphi'}{\vphi}\right)^2 - 2(n-1)\frac{\vphi''}{\vphi}.
\end{equation}
Note that the scalar curvature $R$ of $M$ is given by $R =\vphi' + \rho$. In the first case when $f$ has a unique critical point or $\vphi(0) = 0$, the above equation (\ref{eqnphi2nd}) is equivalently to equation (1.3) in \cite{DSesum} by changing variables. If one chooses the round metric $\bar{g}_{\text{can}}$ on $\sph^{n-1}$ with radius one, i.e., $\vphi'(0) = 1$, then $\bar{R} = (n-1)(n-2)$ and $\rho$ appears as a parameter in the differential equation of $\vphi$. The arguments for their equation (1.3) in \cite{DSesum} which can also be derived from equation (\ref{eqnphi2nd}) show that there is a unique complete Yamabe soliton metric for every $\rho \in \Real$. The asymptotic behavior of $\vphi$ which determines the asymptotic geometries in various cases can also be derived from their arguments.

In this paper we consider the classification of Yamabe solitons in the second case, i.e., the manifold is topologically a product as $M = \Real \times N^{n-1}$ with $(N, \bar{g})$ being of constant scalar curvature $\bar{R}$, and $g = \df r^2 + \vphi^2 \bar{g}$ defines a complete metric on $M$ where $\vphi \in C^{\infty}(\Real)$ is a positive function. We call a Yamabe soliton metric is a \emph{product soliton} if $\vphi$ is a positive constant, i.e., $M$ is isometric to the Riemannian product $\Real \times N$. Note that in \cite{CaoSunZhang} a gradient Yamabe soliton is called \emph{trivial} if the potential function $f$ is constant. Our first result shows that a gradient steady Yamabe soliton is necessarily a product soliton.

\begin{thm}\label{thmsteady}
Any complete gradient steady Yamabe soliton on $M^n = \Real \times_{\vphi} N^{n-1}$ is isometric to the Riemannian product with constant $\vphi$ and $(N, \bar{g})$ being of zero scalar curvature.
\end{thm}

It has the following corollary by using the results above in the case when $f$ has a critical point.
\begin{cor}
Up to scaling, a non-trivial complete gradient steady Yamabe soliton is either the product soliton or the unique one on $\Real^n$ with rotational symmetry.
\end{cor}

On the other hand, by studying the solutions to the differential equation in $\vphi$ we obtained a family of examples for complete gradient Yamabe shrinkers which are not product solitons. More precisely we have
\begin{thm}\label{thmshrinker}
Let $\rho > 0$ be a constant and $n \geq 3$. Suppose $\bar{R} > \rho \frac{n-2}{n+2}\left((n-1)(n+2)\right)^{\frac{4}{n+2}}$ is a positive constant and let $(N^{n-1}, \bar{g})$ be a Riemannian manifold with constant scalar curvature $\bar{R}$. Then there is a unique complete gradient Yamabe soliton metric on the warped product $M^n = \Real \times_{\vphi} N^{n-1}$ with constant $\rho$ that is not isometric to the product metric.
\end{thm}
We briefly discuss the proofs of our results. Since the variable $r$ does not appear explicitly in the equation (\ref{eqnphi2nd}), by introducing new variables the equation can be reduced to a first order nonlinear differential equation. It turns out that the new equation is a classical one, the Abel differential equation of the second kind, see for example \cite{Ince}. This equation does not have an explicit solution in general. So we consider the planar dynamical system defined by this equation. Every trajectory $\sigma$ with positive $\vphi$ of the dynamical system defines a Riemannian metric on $M$. The completeness of the metric requires that certain integral divergences in both directions along $\sigma$. In the steady case we show that all trajectories fail the divergence condition either in one or both directions. This gives us Theorem \ref{thmsteady}. In the shrinking case, for each dimension $n \geq 3$, the integral test of the completeness singles out a unique trajectory $\gamma$ of this system. When $\bar{R}$ and $\rho$ satisfy the equality in Theorem \ref{thmshrinker} instead of the strict inequality, the function $\vphi$ is positive except for one point on $\gamma$. On the other hand when $\bar{R}$ and $\rho$ satisfy the strict inequality, by using some estimates on the trajectory $\gamma$ we can bound $\gamma$ by some simpler curves in the phase plane which allows us to show that $\vphi$ is positive.
 
\begin{rem} 
The equality case of $\bar{R}$ and $\rho$ suggests that $\frac{n-2}{n+2}\left((n-1)(n+2)\right)^{\frac{4}{n+2}}$ is the optimal constant, i.e., we expect that any complete gradient Yamabe shrinker on $M$ is isometric to the Riemannian product with $\vphi$ being a constant function if the inequality of $\bar{R}$ and $\rho$ in Theorem \ref{thmshrinker} does not hold.
\end{rem}

\begin{rem}\label{remRbarpositive}
It was pointed out in \cite{DSesum} that steady and shrinking Yamabe solitons have non-negative scalar curvature. From the second case in the theorem by Cao-Sun-Zhang, we can assume that $\bar{R} > 0$ in the proofs of Theorem \ref{thmsteady} and Theorem \ref{thmshrinker}.
\end{rem}

\begin{rem}
In the expanding case, the Yamabe soliton equation defines dynamical systems that are similar to those in the shrinking and steady cases. The phase-plane analysis of such systems shows that there are complete Yamabe expanders on the warped products with non-constant function $\vphi$. For example when $n=6$, for negative values of $\bar{R}$ there are such examples. The full classification in the expanding case will appear elsewhere.
\end{rem}

The paper is organized as follows. In Section 2 by using new variables the equation (\ref{eqnphi2nd}) in $\vphi$ is reduced to a first order nonlinear differential equation and then we define the dynamical system for this first order equation. In Section 3 we analyze the system for steady Yamabe solitons and prove Theorem \ref{thmsteady}. In Section 4 we study the system for Yamabe shrinkers and prove Theorem \ref{thmshrinker}. In these two sections, one will see that the systems are considerably simpler when $n= 6$. So we prove the results for $n=6$ first and then for other dimensions. There are several figures prepared by the computer algebra system Mathematica which illustrate the ideas of the proofs in this paper. However our arguments do not rely on any of these figures.

\medskip

\textbf{Acknowledgment.} The author would like to thank Huai-Dong Cao for bringing the problems to his attention and many useful suggestions and discussions. He also wants to thank Xiaofeng Sun for helpful conversations.

\medskip

\section{Preliminaries}
In this section we derive the first order differential equations from the equation (\ref{eqnphi2nd}) in $\vphi$ and define a planar dynamical system for this equation. Then the problem of finding a complete gradient Yamabe soliton is equivalent to the one of finding a trajectory satisfies certain restrictions, see Proposition \ref{propDSwzsoliton}.

Since the variable $r$ does not appear in equation (\ref{eqnphi2nd}) explicitly we let $y = \vphi'(r)$ and then we have
\[
\vphi''(r) = \frac{d y}{d \vphi} y = y'(\vphi) y.
\]
Let
\[
x = \vphi(r) > 0, \quad y = \vphi'(r),
\]
and then equation (\ref{eqnphi2nd}) is given by
\[
2(n-1) x y y'(x) = \bar{R} - (n-1)(n-2)y^2 - x^2(y + \rho).
\]
It can be rewritten as
\[
y y'(x) = - \frac{n-2}{2x}y^2 - \frac{x}{2(n-1)}y + \frac{\bar{R}}{2(n-1)x} - \frac{\rho x}{2(n-1)}.
\]
The above equation has the form of the Abel differential equation of the second kind
\[
y y'(x) = a_1(x) y^2 + a_2(x) y + a_3(x)
\]
with
\[
a_1 (x) = - \frac{n-2}{2x}, \quad a_2(x) = - \frac{x}{2(n-1)}, \quad a_3(x) = \frac{\bar{R}}{2(n-1)x} - \frac{\rho x}{2(n-1)}.
\]
Using the standard process we write the equation in the canonical form. Let
\[
y = E(x) w, \quad \mbox{where}\quad E(x) = \exp\left(\int a_1(x) dx \right) = x^{- \frac{n-2}{2}},
\]
then we have
\[
w w'(x) = F_1(x) w + F_0(x)
\]
with
\[
F_1(x) = \frac{a_2(x)}{E(x)} = - \frac{1}{2(n-1)}x^{\frac{n}{2}}, \quad F_0(x) = \frac{a_3(x)}{E^2(x)} = \frac{\bar{R}}{2(n-1)}x^{n-3} - \frac{\rho}{2(n-1)}x^{n-1}.
\]
By introducing the new variable
\[
z = - \int F_1(x) dx = \frac{1}{(n-1)(n+2)}x^{\frac{n+2}{2}}
\]
the equation for $w$ has the following canonical form
\begin{equation*}
w w'(z) + w = \Phi(z)
\end{equation*}
where
\begin{eqnarray*}
\Phi(z) & = & - \frac{F_0(x)}{F_1(x)} = \bar{R} x^{\frac{n}{2}-3} - \rho x^{\frac{n}{2}-1} \\
& = & \bar{R} \left((n-1)(n+2)\right)^{\frac{n-6}{n+2}}z^{\frac{n-6}{n+2}} - \rho \left((n-1)(n+2)\right)^{\frac{n-2}{n+2}}z^{\frac{n-2}{n+2}}.
\end{eqnarray*}
The transformation formulas between $(x, y)$ and $(z, w)$ are given by
\begin{equation}\label{eqnTrfmzwxy}
\left\{
\begin{array}{rcl}
z & = & \dfrac{1}{(n-1)(n+2)}x^{\frac{n+2}{2}} \\
w & = & y x ^{\frac{n-2}{2}},
\end{array}
\right.
\end{equation}
and
\begin{equation}\label{eqnTrfmxyzw}
\left\{
\begin{array}{rcl}
x & = & \left((n-1)(n+2)z\right)^{\frac{2}{n+2}} \\
y & = & w \left((n-1)(n+2)z\right)^{- \frac{n-2}{n+2}}.
\end{array}
\right.
\end{equation}

We summarize the above discussion as
\begin{prop}\label{propDSwzsoliton}
A gradient Yamabe soliton metric $g$ on $M$ is determined uniquely by one of the followings:
\begin{enumerate}
\item A positive solutions $\vphi = \vphi(r)$ to equation (\ref{eqnphi2nd}).
\item A solution $w = w(z)$ to the following equation defined on positive $z$-axis:
\begin{equation}\label{eqnwzPhi}
w w'(z) + w = \Phi(z)
\end{equation}
with
\[
\Phi(z) = \bar{R}\left((n-1)(n+2)z\right)^{\frac{n-6}{n+2}} - \rho \left((n-1)(n+2) z \right)^{\frac{n-2}{n+2}}.
\]
\item A trajectory $\sigma =(z(s), w(s))$ on the half plane with positive $z$ of the following dynamical system:
\begin{equation}\label{eqnDSzw}
\left\{
\begin{array}{rcl}
\dot{z}(s) & = & w \\
\dot{w}(s) & = & \Phi(z) - w,
\end{array}
\right.
\end{equation}
where $\Phi(z)$ is given in (2) above.
\end{enumerate}
Moreover the metric $g$ is complete if and only if either $\vphi= \vphi(r)$ is defined for all $r \in \Real$, or the following integral
\begin{equation}\label{eqnIcompleteness}
I = \int_{\sigma} \frac{1}{w}z^{- \frac{2}{n+2}}\df z
\end{equation}
diverge in both directions along a trajectory $\sigma$ in (3).
\end{prop}
\begin{proof}
We already showed the equivalence of (1) and (2). The dynamical system formulation of (2) gives us (3). Next we consider the completeness of the metric $g$ defined by $\vphi$ or a trajectory $\sigma$ (3). The metric $g$ is complete if and only if its restriction on the real line $\Real$ which is a geodesic in $M$ is complete, i.e.,  $\vphi(r)$ exists for all $r\in (- \infty,\infty)$ and stays positive. Suppose $\sigma = (z(s), w(s))$ is defined for $s \in (s_0, s_1)$ where $s_0$ or $s_1$ may be infinity and let $(z_i, w_i) = \lim_{s\rightarrow s_i}(z(s), w(s))$ for $i=0, 1$. Take another point $(a,b)$ on $\sigma$ between $(z_0, w_0)$ and $(z_1, w_1)$. Since $\df r = \frac{1}{y} \df x$ which is a constant multiple of $w^{-1}z^{-\frac{2}{n+2}}\df z$, the completeness of $g$ is equivalent to the divergence of the integral $I$ in both directions along $\sigma$. Note that it diverges to infinity with different signs.
\end{proof}

\begin{rem}
In the rest of the paper, if two integrals $I_1$ and $I_2$ diverge or converge simultaneously along a trajectory, then we write $I_1 \sim I_2$.
\end{rem}

\begin{rem}
In the dynamical system (\ref{eqnDSzw}) since $\df s = \frac{1}{w}\df z$ and $\df r = w^{-1}z^{-\frac{2}{n+2}}\df z$, we have
\begin{equation}\label{eqndsdr}
\df s = z^{\frac{2}{n+2}}\df r,
\end{equation}
i.e., $s$ points in the positive direction of $r$ and it is comparable with $r$ when the trajectory approaches to a point $(z_0, w_0)$ with $z_0 > 0$.
\end{rem}

\begin{rem}
If we use the variables $x$ and $y$ instead of $z$ and $w$, then the dynamical system is given by
\begin{equation}\label{eqnDSxy}
\left\{
\begin{array}{rcl}
\dot{x}(t) & = & 2(n-1)x y \\
\dot{y}(t) & = & \bar{R} - (n-1)(n-2)y^2 - x^2 (y + \rho),
\end{array}\right.
\end{equation}
and the integral for the completeness of $g$ is
\[
\int_{\sigma} \frac{1}{y}\df x
\]
for a trajectory $\sigma$ in $xy$-plane with positive $x$.
\end{rem}

\begin{rem}
The system (\ref{eqnDSxy}) is more convenient in the case when $\vphi$ has a zero point, i.e., the metric is rotationally symmetric. We give a brief argument showing that there is a unique complete rotationally symmetric metric $g$ in this case based on the phase plane analysis. The initial condition that $\vphi(0) = 0$ and $\vphi'(0) = 1$ indicates that we should look at the trajectories that approach to the point $(0,1)$ in the $xy$-plane. This point is a critical point as $\dot{x} = \dot{y} = 0$. The linearization of the system at this point is given by
\[
\begin{pmatrix}
\dot{x} \\
\dot{y}
\end{pmatrix}
=
\begin{pmatrix}
2(n-1) & 0 \\
0 & -2(n-1)(n-2)
\end{pmatrix}
\begin{pmatrix}
x \\
y
\end{pmatrix}.
\]
It follows that $(0,1)$ is a saddle point. One can see that two separatrices are the $y$-axis, and the separatrix $\gamma$ with horizontal tangent and positive $x$ defines the Riemannian metric $g$. The trajectory $\gamma$ approaches to another critical point as $t\rightarrow \infty$ if $\rho > 0$ and to infinity if $\rho \leq 0$. In either case, one can show that the integral $\int y^{-1}\df x$ diverges as $t\rightarrow \infty$ which gives the completeness of $g$.
\end{rem}

\medskip

\section{Gradient Steady Yamabe solitons}
In this section we prove Theorem \ref{thmsteady}, i.e., any complete gradient steady Yamabe soliton on a warped product is necessarily a Riemannian product with constant $\vphi$. From the theorem by Cao-Sun-Zhang and Remark \ref{remRbarpositive} in the introduction, we only have to show

\begin{thm}\label{thmsteadyRbarpos}
Let $n \geq 3$ and $(N^{n-1}, \bar{g})$ be a Riemannian manifold with positive constant scalar curvature. Then there is no complete gradient steady Yamabe soliton metric on $M^n = \Real\times_{\vphi} N^{n-1}$.
\end{thm}
Our proof of the theorem above is separated into two cases, the case with $\dim M = 6$, see Theorem \ref{thmsteadyN6}, and the case with other dimensions, see Theorem \ref{thmsteadyNgeneral}.

First we fix the notions and state some general facts for all dimensions. Since $M$ is steady soliton, we have $\rho = 0$ and thus
\[
\Phi(z) =\lambda z^{\frac{n-6}{n+2}} \quad \text{where} \quad \lambda =  \left((n-1)(n+2)\right)^{\frac{n-6}{n+2}}\bar{R} > 0.
\]
By re-scaling the metric $g$ if necessary we may assume that $\lambda = 1$. So equation (\ref{eqnwzPhi}) and the dynamical system (\ref{eqnDSzw}) have the following form:
\begin{equation}\label{eqnwzsteady}
w(z) w'(z) + w(z) = z^{\frac{n-6}{n+2}},
\end{equation}
and
\begin{equation}\label{eqnDSzwsteady}
\left\{
\begin{array}{rcl}
\dot{z}(s) & = & w \\
\dot{w}(s) & = & z^{\frac{n-6}{n+2}} - w.
\end{array}\right.
\end{equation}

Let
\[
\mathcal{P} = \set{(z,w) \in \Real^2 : z > 0}
\]
denote the half plane where we are looking for trajectories. We define the following curve $S$:
\[
S = \set{(z, w) \in \mathcal{P} : w = \Phi(z)}.
\]
It separates the half plane into two pieces:
\[
\mathcal{P}_1 = \set{(z,w) \in \mathcal{P} : w > \Phi(z)}, \quad \mathcal{P}_2 = \set{(z, w) \in \mathcal{P} : w < \Phi(z)}.
\]

In any dimension we show that a trajectory does not define a complete Riemannian metric when it approaches the $w$-axis, the boundary of $\mathcal{P}$.

\begin{lem}\label{lemsteadywaxis}
Suppose $\sigma$ is a trajectory of the system (\ref{eqnDSzwsteady}) that approaches a point $p_0 = (0, w_0)$ on the $w$-axis. Then the integral $I$ along $\sigma$ does not diverge near $p_0$.
\end{lem}
\begin{proof}
If $p_0$ is not the origin, i.e., $w_0 \ne 0$, then in the integrand of $I$, $\abs{w^{-1}}$ is bounded from infinity near the point $p_0$. Since $\frac{2}{n+2} < 1$, the integral $I$ converges. If $p_0$ is the origin, in the case when $n= 6$ we can use the explicit solution to the equation (\ref{eqnwzsteady}) as $\Phi(z) = 1$, see for example \cite{PolyaninZaitsev}, to see that $I$ converges near the origin. In the following we give a general argument for all dimensions.

Let $u = z^{\frac{1}{n+2}}$ as $z\geq 0$ and then the equation (\ref{eqnwzsteady}) of $z$ and $w$ is equivalent to
\[
\frac{\df w}{\df u} = \frac{n+2}{w}\left(u^{2n-5} - u^{n+1} w \right),
\]
and the integral $I$ is given by
\[
I = (n+2)\int\frac{u^{n-1}}{w} \df u.
\]
The above equation of $u$ and $w$ defines the following dynamical system
\begin{equation}\label{eqnDSuw}
\left\{
\begin{array}{rcl}
\dot{u}(t) & = & w \\
\dot{w}(t) & = & (n+2)\left(u^{2n-5} - u^{n+1} w\right)
\end{array}\right.
\end{equation}

If $n = 3$, then the linearization of system (\ref{eqnDSuw}) has the following form
\[
\begin{pmatrix}
\dot{u}(t) \\
\dot{w}(t)
\end{pmatrix}
=
\begin{pmatrix}
0 & 1 \\
5 & 0
\end{pmatrix}
\begin{pmatrix}
u \\
w
\end{pmatrix}.
\]
The above coefficient matrix has eigenvalues $\pm\sqrt{5}$ with eigenvectors $(\pm 1, \sqrt{5})$. So the origin is the saddle point and there two different trajectories approach it, one with $t\To -\infty$ and the other with $t \To \infty$. We use the power series to approximate these two trajectories. Let
\[
w = u^k(a_0 + a_1 u + \cdots )
\]
with $k> 0$ and $a_0 \ne 0$. Then comparing the both sides of the equation on $u$ and $w$ shows that $k = 1$ and $a_0^2 = 5$. So for some $a> 0$ small enough, we can estimate the integral $I$ as
\[
I = (n+2)\int_0^{a} \frac{u^2}{w}\df u \sim \int_0^a \frac{u^2}{u}\df u = \int_0^a u \df u < \infty.
\]

If $n\geq 4$ then the linearization of the system (\ref{eqnDSuw}) has the coefficient matrix
\[
\begin{pmatrix}
0 & 1 \\
0 & 0
\end{pmatrix}
\]
which shows that $(0,0)$ is a non-hyperbolic critical point of the nonlinear system. However since the coefficient matrix is not a zero matrix, by \cite[pp. 180, Theorem 2]{Perko} we know that $(0,0)$ is a topological saddle, i.e., there are two different trajectories in $\mathcal{P}$ approaching this point. A power series approximation as in the case when $n=3$ shows that $w = u^{n-2}(a_0 + a_1 u + \cdots)$ with nonzero $a_0$ around the origin and then $I$ is finite when $z\To 0$.
\end{proof}

\subsection{The case with $\dim M = 6$} We have $\Phi(z) = 1$ and $S$ is horizontal line $w = 1$.
The following Figure \ref{figpsteadyn6} shows the phase portrait of the dynamical system (\ref{eqnDSzwsteady}). Note that the horizontal axis is the $z$-axis and the vertical one is the $w$-axis. The green line is $S$ and the half plane $\mathcal{P}$ is separated into two regions.

\begin{figure}[!htb]
\begin{center}
\includegraphics[width=200pt,height=200pt]{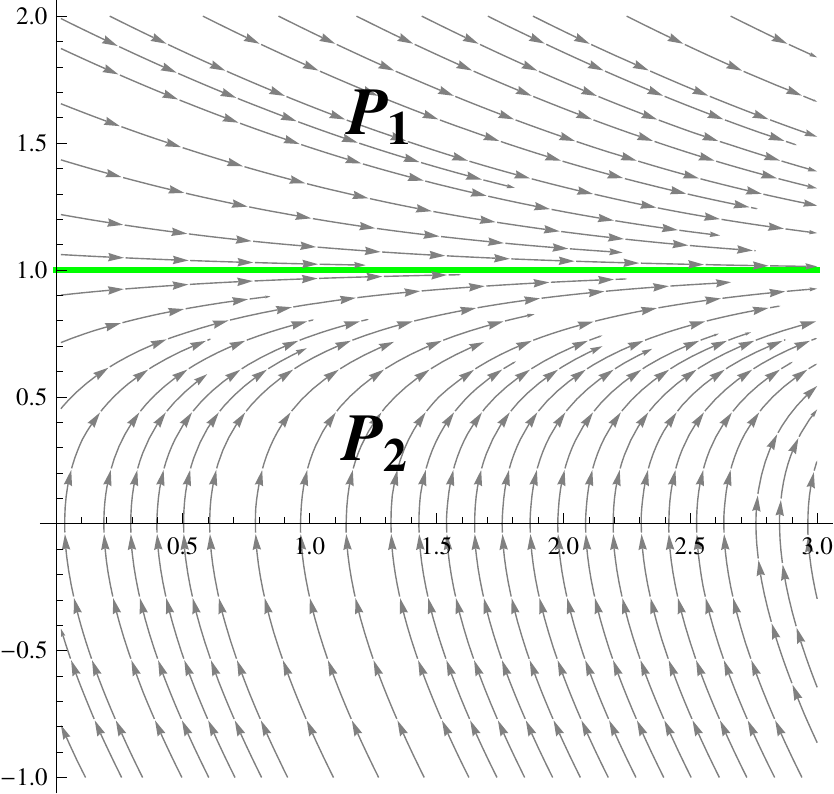}
\end{center}
\caption{\footnotesize Phase portrait with $n=6$, $\lambda=1$ and $\rho = 0$}\label{figpsteadyn6}
\end{figure}

Note that $S$ is a trajectory of the system (\ref{eqnDSzwsteady}) and then any trajectory passing through a point in $\mathcal{P}_1$(or $\mathcal{P}_2$) will stay in that region.

\begin{thm}\label{thmsteadyN6}
There is no complete steady Yamabe soliton on $M^6 = \Real \times_{\vphi} N^5$ with $N^5$ being positive scalar curvature.
\end{thm}

\begin{proof}
To show the statement we only have to show that any trajectory cannot define a complete metric. Note that there is no vertical asymptote in the half plane $\mathcal{P}$. Since along a trajectory we have $w'(z) = \frac{1}{w} -1$, the $w$-axis is not a vertical asymptote of the trajectory. First the trajectory $w=1$ does not define a complete metric. Next we consider the trajectory that lies in the region $\mathcal{P}_1$. Suppose $\gamma$ is a trajectory passing through a point $\gamma(0) = (z_0, w_0) \in \mathcal{P}_1$. Then $w(z)$ is a decreasing function in $z$ and it meets $w$-axis with positive value of $w$. This shows that the metric is not complete as $z \To 0$.

Now suppose $(z_0, w_0) \in \mathcal{P}_2$ and $w_0 > 0$. Since $w$ is an increasing function in $z$, $\gamma$ meets either $w$-axis or $z$-axis as $z$ decreases. If it meets $w$-axis(including the origin), then Lemma \ref{lemsteadywaxis} shows that the metric is not complete when $z\To 0$. If $\gamma$ meets $z$-axis with positive value, then it intersects with $z$-axis vertically and thus enters the 4th quadrant. In this quadrant $w$ is a decreasing function in $z$. Since $w'(z) = \frac{1-w}{w} < -1$ we have $w(z) < - z$ for $z$ large enough, i.e., $- w > z$. For the integral $I$ when $z\To \infty$ we have
\[
I = \int_{z_0}^{\infty} \frac{1}{w}z^{-\frac{1}{4}}\df z > \int_{z_0}^{\infty} - \frac{1}{z}z^{-\frac{1}{4}}\df z = - \int_{z_0}^{\infty} z^{-\frac{5}{4}}\df z > -\infty.
\]
So the trajectory does not give a complete metric. This finishes the proof.
\end{proof}

\subsection{The other case with $\dim M \ne 6$} 
\begin{figure}[!htb]
\begin{center}
\includegraphics[width=400pt,height=200pt]{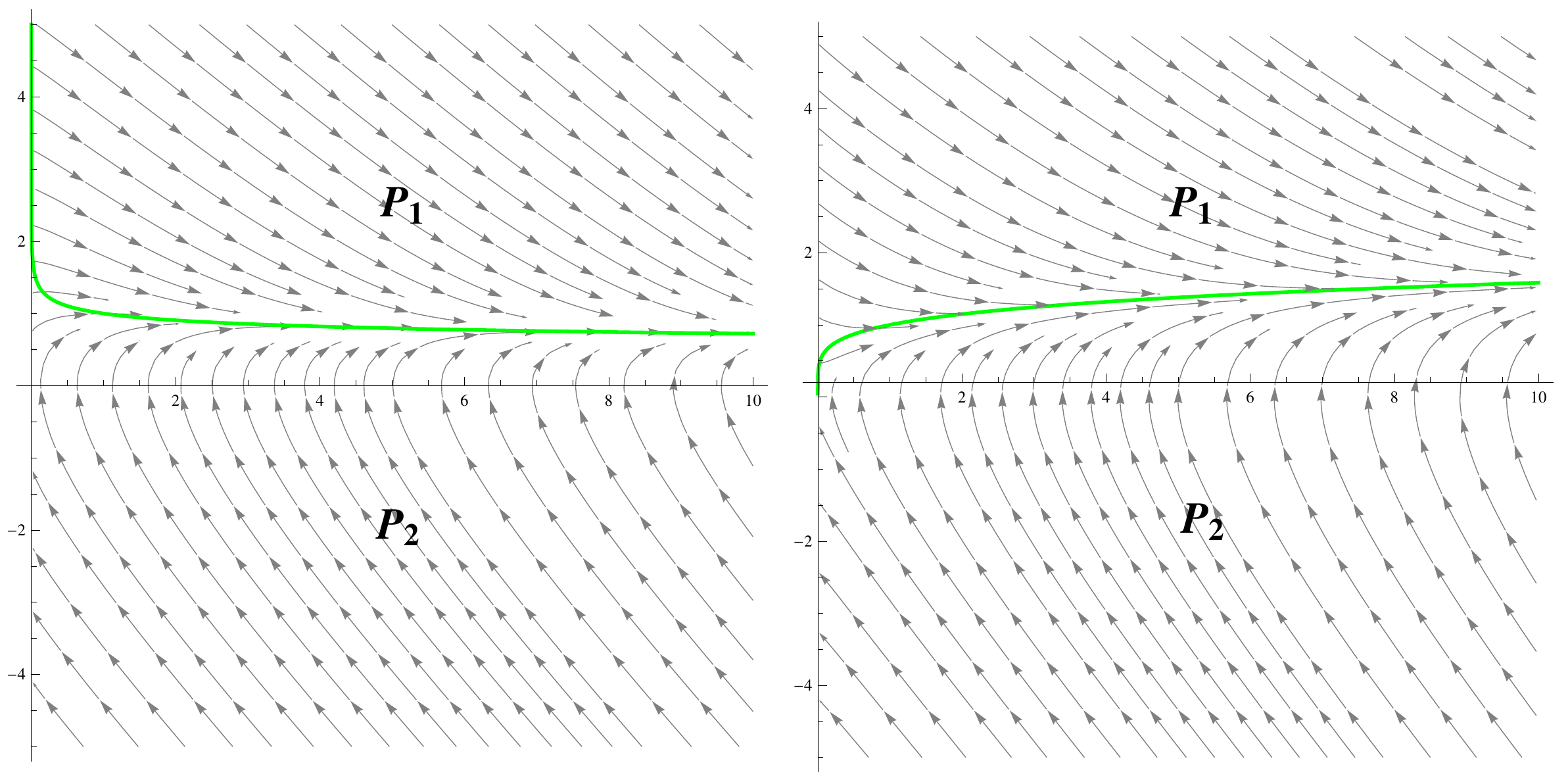}
\end{center}
\caption{\footnotesize Phase portraits with $n=5$(the left one) and $n=8$(the right one), $\lambda=1$ and $\rho = 0$}\label{figpsteadyN58}
\end{figure}

We have $\Phi(z) = z^{\frac{n-6}{n+2}}$ and it is zero or approaches $\infty$ as $z\To 0$ when $n > 6$ or $n< 6$. The typical phase portraits in this case are shown in Figure \ref{figpsteadyN58} for $n=5$(the one on the left) and $n=8$(the one on the right). The green curve is $S$ for each case.

\begin{thm}\label{thmsteadyNgeneral}
For any $n \ne 6$ there is no complete steady Yamabe soliton on $M^n = \Real \times_{\vphi} N^{n-1}$ with $N^{n-1}$ being positive scalar curvature.
\end{thm}
\begin{proof}
The proof is similar to the one when $n=6$. Since $S$ is not a trajectory, any trajectory starts from a point on $S$ will enter the region $\mathcal{P}_1$ or $\mathcal{P}_2$. Note that there is no vertical asymptote of trajectories in the half plane $\mathcal{P}$. When $n > 6$, the $w$-axis is not a vertical asymptote either. We claim that the $w$-axis is not a vertical asymptote when $n < 6$. Suppose not, there is a trajectory $\gamma = (z, w(z))$ in the region $\mathcal{P}_1$ such that along $\gamma$ we have
\[
\lim_{z\To 0} w = \infty.
\]
Since in the region $\mathcal{P}_2$ with positive $w$, $w$ decreases as $z$ tends to zero, $\gamma$ is bounded from below by the curve $S$. In particular it follows that there exists a sequence $\set{z_i > 0}$ that converges to zero, such that $\Phi'(z_i) \geq w'(z_i)$, i.e.,
\[
\frac{n-6}{n+2}z_i^{- \frac{8}{n+2}} \geq \frac{z_i^{\frac{n-6}{n+2}}}{w(z_i)} - 1.
\]
Multiplying $z_i^{\frac{8}{n+2}}$ on both sides yields
\[
\frac{n-6}{n+2} \geq \frac{1}{w(z_i)} - z_i^{\frac{8}{n+2}}.
\]
Let $i \To \infty$, then the right hand side converges to zero. However the left hand side is a fixed negative number which shows a contradiction.

As in the proof of case $n= 6$, from the monotonicity properties of the function $w = w(z)$ in different regions, one can show that a trajectory either meets the $w$-axis with finite $w$, or enters the region $\mathcal{P}_2$ with negative $w$ and then decreases at least like the function $w = -z$ for $z$ large. In either case, the integral $I$ is finite and so the metric defined by the trajectory is not complete which finishes the proof.
\end{proof}

\medskip
\section{Gradient Yamabe Shrinkers}

In this section we prove Theorem \ref{thmshrinker}, the existence of complete gradient Yamabe shrinkers on warped products which are not isometric to the Riemannian product. We consider the case $n=6$ first and prove Theorem \ref{thmshrinkerN6}, and then the case with other dimensions and prove Theorem \ref{thmshrinkerNgeneral}.

First we fix notions for this section and show some general facts for all dimensions. Let $\bar{R} > 0$ be the scalar curvature of $N^{n-1}$, $A = \left((n-1)(n+2)\right)^{\frac{n-6}{n+2}}$ and $\lambda = \bar{R}A > 0$. We may assume that $\rho = \left((n-1)(n+2)\right)^{-\frac{n-2}{n+2}} > 0$ by re-scaling if necessary and then we have
\[
\Phi(z) = \lambda z^{\frac{n-6}{n+2}} - z^{\frac{n-2}{n+2}}.
\]
So the differential equation (\ref{eqnwzPhi}) in $z$ and $w$ is given by
\begin{equation}\label{eqnwzshrinker}
w(z)w'(z) + w(z) = \lambda z^{\frac{n-6}{n+2}} - z^{\frac{n-2}{n+2}},
\end{equation}
and the dynamical system (\ref{eqnDSzw}) has the following form
\begin{equation}\label{eqnDSzwshrinker}
\left\{
\begin{array}{rcl}
\dot{z}(s) & = & w \\
\dot{w}(s) & = & \lambda z^{\frac{n-6}{n+2}} - z^{\frac{n-2}{n+2}} - w.
\end{array}\right.
\end{equation}

First we observe that for all $\lambda$ and $n$, $(\xi, 0)$ is a critical point of the dynamical system (\ref{eqnDSzwshrinker}) where
\[
\xi = \lambda^{\frac{n+2}{4}}.
\]
If $n \geq 7$, then $(0,0)$ is another critical point of the system. In the case of the steady solitons, we showed that the integral for completeness test is finite along a trajectory when it approaches the origin. Similarly we have

\begin{lem}
Suppose $\gamma$ is a trajectory in $\mathcal{P}$ that approaches the origin $(0,0)$, then the integral $I$ is finite.
\end{lem}
\begin{proof}
We use the variable $u = z^{\frac{1}{n+2}}$ around the origin $(0,0)$. The differential equation (\ref{eqnwzshrinker}) of $z$ and $w$ can be written as
\[
\frac{\df w}{\df u} = (n+2)\frac{\lambda u^{2n-5} - u^{2n-1}}{w} - (n+2) u^{n+1}.
\]
Following the argument of Lemma \ref{lemsteadywaxis} in the steady case, we know that $(0,0)$ is a topological saddle of the dynamical system in $u$ and $w$.
A formal power series expansion shows that
\[
w = u^{n-2}(a_0 + a_1 u + \cdots)
\]
where $(n-2)a_0^2 = (n+2)\lambda$. So for some $a> 0$ small we have
\[
I = \int_0^a \frac{1}{w z^{\frac{2}{n+2}}}\df z \sim \int_{0}^{a}\frac{1}{z^{\frac{n}{n+2}}}\df z < \infty
\]
which finishes the proof.
\end{proof}

The next two propositions characterize the local behavior of the trajectories around the critical point $(\xi,0)$.

\begin{prop}
The nonlinear system (\ref{eqnDSzwshrinker}) has $(\xi, 0)$ as either a stable node when $(n+2)\lambda \geq 16$, or a sable focus when $(n+2)\lambda < 16$. In both cases the integral $I$ diverges to infinite along a trajectory that approaches this point.
\end{prop}
\begin{proof}
At the point $(\xi,0)$, the nonlinear system (\ref{eqnDSzwshrinker}) has the following linearization
\begin{equation}\label{eqnDSzwlinear}
\begin{pmatrix}
\dot{z} \\
\dot{w}
\end{pmatrix}
=
\begin{pmatrix}
0 & 1 \\
-\frac{4}{(n+2)\lambda } & -1
\end{pmatrix}
\begin{pmatrix}
z \\
w
\end{pmatrix}.
\end{equation}
Let $r_1$ and $r_2$ be eigenvalues of the coefficient matrix of the above linear system. We have
\[
r_1 + r_2 = -1\quad \text{and} \quad r_1 r_2 = \frac{4}{(n+2)\lambda}.
\]
If $(n+2)\lambda - 16 \geq 0$ then the two eigenvalues are negative real numbers and so $(\xi, 0)$ is a stable node of the linear system. If $(n+2)\lambda -16 < 0$, then the two eigenvalues are complex with negative real part and so $(\xi, 0)$ is a stable focus. On the right hand side of the nonlinear system (\ref{eqnDSzwshrinker}), we only have power functions in $z$ and $w$ with positive exponents, from \cite[pp. 142, Theorem 4]{Perko} the critical point $(\xi,0)$ is a stable node or a stable focus of the nonlinear system respectively. In both cases if $\gamma$ is a trajectory around this point, then $\gamma$ approach $(\xi, 0)$ as $s \To \infty$. Since the distance function $r$ is comparable with the parameter $s$ near the critical point, $r$ tends to infinite as $\gamma$ approaches this point.
\end{proof}

In general we only know that a trajectory near the point $(\xi, 0)$ will approach it as $s\To \infty$. In this case we can show the convergence for a large area.

\begin{prop}\label{propshrinkercriticalptxi}
For any $0 < z_0 < \xi$, the trajectory with initial point $(z_0, 0)$ converges to the critical point $(\xi,0)$ as $s\To \infty$.
\end{prop}
\begin{proof}
Suppose $\gamma$ is a such trajectory with $\gamma(0) = (z_0, 0)$. As $s$ becomes positive, $\gamma$ enters the first quadrant and then $w=w(z)$ is an increasing function of $z$. As $\dot{z}(s) = w > 0$, $z$ and $w$ are increasing as $s$ increases. Then $\gamma$ will meet $S_1$ and after that $z$ increases while $w$ decreases. The trajectory $\gamma$ will meet the $z$-axis and enter the 4th quadrant. Then it will meet $S_1$ and the $z$-axis again. Suppose $\gamma$ meets the $z$-axis from the 4th quadrant at $s=s_1$, i.e., $\gamma(s_1) = (z_b,0)$.

We claim that $z_b > z_0$. Suppose not, since there is no periodic orbit by Bendixson's criteria, see for example \cite[pp. 245, Theorem 1]{Perko}, we have $z_b < z_0$. Then the curve $\gamma(s)$ with $s \in [0, s_1]$ and the segment $\set{(z, 0) : z_b \leq z \leq z_0}$ form a bounded region and $\gamma(s)$ stays in this region when $s< 0$. So $\gamma(s)$ will approach a critical point as $s\To -\infty$. This contradicts the fact that in this bounded region there is only one critical point $(\xi, 0)$ which is stable.

Now we have $z_b > z_0$ and we have a bounded region $E$ whose boundary consists of the curve $\gamma(s)$ with $s\in [0,s_1]$ and the segment $\set{(z, 0): z_0 \leq z \leq z_b}$. The trajectory $\gamma(s)$ for $s\geq s_1$ stays in $E$ and so $\gamma(s)$ approaches the unique critical point $(\xi,0)$ as $s\To \infty$.
\end{proof}

In the following we show that there is a unique trajectory $\gamma(s) = (z(s), w(s))$ that approaches $(\xi, 0)$ as $s \To \infty$ and the function $w = w(z)$ has a distinguished asymptotic expansion at infinity when $s \To -\infty$. Such asymptotic expansion of $w$ implies that the integral $I$ diverges when $s \To -\infty$ and so the metric defined by $\gamma$ is complete.

\subsection{The case with $\dim M =6$} We have
\[
\Phi(z) = \lambda - \sqrt{z},
\]
and
\[
S_1 = \set{(z,w) \in \mathcal{P}: w = \lambda - \sqrt{z}}.
\]
The curve $S_1$ determines the monotonicity of trajectories. We introduce another curve in the half plane $\mathcal{P}$ which determines the convexity of trajectories. From the formula $w''(z)$ of a trajectory $\sigma = (z,w(z))$ we define
\[
S_2 = \set{(z,w)\in \mathcal{P} : \frac{w^2}{\sqrt{z}} - 2(\lambda - \sqrt{z})w + 2(\lambda - \sqrt{z})^2 = 0 }.
\]

The curve $S_1$ separates the phase plane $\mathcal{P}$ into two pieces:
\[
\mathcal{P}_1 = \set{(z,w)\in \mathcal{P} : w > \lambda -\sqrt{z}}, \quad \mathcal{P}_2 = \set{(z,w) \in \mathcal{P} : w < \lambda -\sqrt{z}}.
\]
Note that $S_1$ is not a trajectory of the system in $z, w$. Any trajectory starts from $S_1$ will enter one of these two regions.

 Let
\begin{eqnarray*}
f_1 & = & (\lambda - \sqrt{z})\sqrt{z} + (\sqrt{z} - \lambda)\sqrt{z - 2\sqrt{z}} \\
f_2 & = & (\lambda - \sqrt{z})\sqrt{z} - (\sqrt{z} - \lambda)\sqrt{z - 2\sqrt{z}}
\end{eqnarray*}
for $z \geq \max\set{4, \lambda^2}$. On the 4th quadrant, $S_2$ has two components given by different formulas:
\begin{eqnarray*}
S_{2a} & = & \set{(z, w)\in \mathcal{P} : w = f_1, \, z \geq \max(4, \lambda^2)} \\
S_{2b} & = & \set{(z, w)\in \mathcal{P} : w = f_2, \, z \geq \max(4, \lambda^2)}.
\end{eqnarray*}
These two curves have a common point $T$ when $z = \max\set{4, \lambda^2}$. If $\lambda \geq 0$, $T$ is the same as the critical point $(\xi, 0)$. If $\lambda < 2$, then $T$ has the coordinate $(4,2\lambda -4)$. It is easy to see that when $z$ is large, the defining function of $S_{2a}$ is asymptotic to $w=-\sqrt{z}$ and the one for $S_{2b}$ is asymptotic to $w = -2 z$. We define the following trapping region $\mathcal{T}$ in $\mathcal{P}_2$:
\[
\mathcal{T} = \set{(z, w) \in \mathcal{P}_2 : f_2 < w < f_1 }.
\]

The following Figure \ref{figRegionTShrinkern6} shows the region $\mathcal{T}$ for some typical values of $\lambda$. The green curve is $S_1$ and the dotted blue curve is $S_2$.
\begin{figure}[!htb]
\begin{center}
\includegraphics[width=480pt,height=160pt]{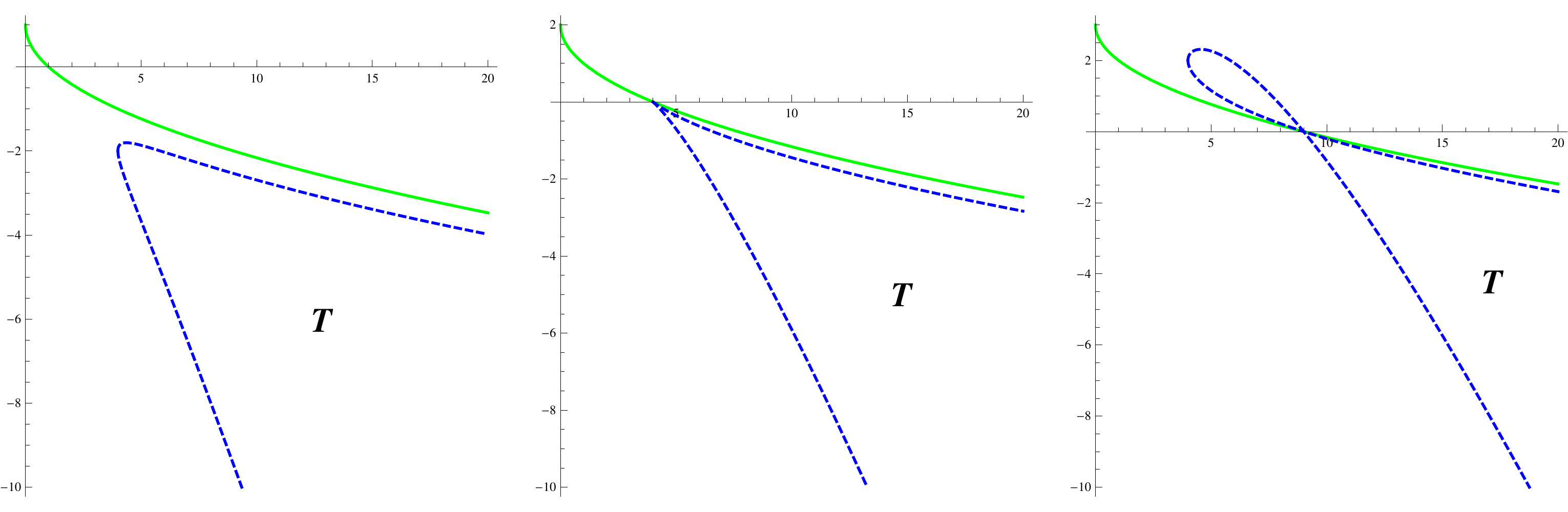}
\end{center}
\caption{\footnotesize The region $\mathcal{T}$ for $\lambda = 1, 2$ and $3$ when $n = 6$.}\label{figRegionTShrinkern6}
\end{figure}

\begin{prop}
For these two curves $S_1$ and $S_2$ we have
\begin{enumerate}
\item If $\lambda \leq 2$, then on $S_2$ we have $z \geq 4$ and the whole curve $S_2$ is in the 4th quadrant as $S_2 = S_{2a}\cup S_{2b}$ and it is below $S_1$.
\item If $\lambda > 2$ then $S_2$ has non-empty components in the 1st quadrant for $4 \leq z < \lambda^2$, and in the 4th quadrant for $z > \lambda^2$ which is $S_{2a}\cup S_{2b}$. The curve $S_2$ is above $S_1$ in the 1st quadrant, and is below $S_1$ in the 4th quadrant.
\end{enumerate}
Furthermore if a trajectory meets the boundary of $\mathcal{T}$ with intersection point $(z_0, w_0)$, then it enters and stays inside the region $\mathcal{T}$ for $z > z_0$.
\end{prop}

\begin{proof}
The defining equation of $S_2$ is given by $w''(z) = 0$ and we can rewrite it as
\[
w''(z) = \frac{1}{-\sqrt{z} w^3}\left(\left(w-(\lambda - \sqrt{z})\sqrt{z}\right)^2 - \sqrt{z}(\sqrt{z}-2)(\lambda - \sqrt{z})^2\right).
\]
It is easy to see that $w=0$ implies that $z = \lambda^2$ for $z> 0$. So the set of $w''(z) = 0$ with positive $z$ consists of the single point $(\lambda^2, 0)$ and the curve $S_2$. At any point on $S_2$ we have
\[
\sqrt{z}(\sqrt{z} - 2)(\lambda - \sqrt{z})^2 = \left(w - (\lambda - \sqrt{z})\sqrt{z}\right)^2 \geq 0.
\]
It follows that $\sqrt{z} -2 \geq 0$, i.e., $z\geq 4$. Moreover the solution to the equation $w''(z) = 0$ is given by
\[
w_{1,2} = (\lambda - \sqrt{z})\sqrt{z} \pm \abs{\sqrt{z}-\lambda}\sqrt{z - 2 \sqrt{z}}.
\]
If $\lambda \leq 2$, then $\sqrt{z}-\lambda \geq 0$ and
\[
w_{1,2} = (\sqrt{z} - \lambda)\left(-\sqrt{z} \pm \sqrt{z- 2\sqrt{z}}\right) \leq 0,
\]
i.e., $S_2$ lies in the 4th quadrant. If $\lambda > 2$, then we have
\[
w_{1,2} = (\lambda - \sqrt{z})\left(\sqrt{z} \pm \sqrt{z- 2\sqrt{z}}\right) \geq 0 \quad \text{if } 4 \leq z < \lambda^2,
\]
and
\[
w_{1,2} = (\sqrt{z}- \lambda )\left(- \sqrt{z} \pm \sqrt{z- 2\sqrt{z}}\right) \leq 0 \quad \text{if } z \geq \lambda^2.
\]
So the curve $S_2$ has nonempty components in the 1st and 4th quadrants.
Using the formulas of $w_{1,2}$ the relative position between $S_1$ and $S_2$ can be checked by the sign of $\lambda - \sqrt{z} - w_{1,2}$ when $z> \lambda^2$(and $4 < z < \lambda^2$ if $\lambda > 2$).

Next suppose $\sigma$ is a trajectory that touch $S_{2a}$ at the point $p=(z, w)$. Since $\dot{z}(s) = w < 0$ in the 4th quadrant, $z$ increases as $s$ decreases. Let $m_1=w'(z)$ be the slope of the tangent of $\sigma$ at $p$ and $m_2 = f_1'(z)$ be the slope for the curve $S_{2a}$. We have
\begin{eqnarray*}
m_1 & = & \frac{\lambda - \sqrt{z}}{f_1} - 1 = \frac{\lambda - \sqrt{z}}{(\lambda - \sqrt{z})\left(\sqrt{z} - \sqrt{z - 2\sqrt{z}}\right)} - 1 = \frac{1}{\sqrt{z} - \sqrt{z - 2\sqrt{z}}} - 1 \\
& = & \frac{1}{2} + \frac{\sqrt{z-2\sqrt{z}}}{2\sqrt{z}} -1, \\
m_2 & = & -1 + \frac{\lambda}{2\sqrt{z}} + \frac{\sqrt{z - 2\sqrt{z}}}{2\sqrt{z}} + \frac{\sqrt{z}-\lambda}{2\sqrt{z - 2\sqrt{z}}}\left(1-\frac{1}{\sqrt{z}}\right),
\end{eqnarray*}
and then
\begin{eqnarray*}
m_2 - m_1 & = & \frac{\lambda}{2\sqrt{z}} - \frac{1}{2} + \frac{\sqrt{z}-\lambda}{2\sqrt{z-2\sqrt{z}}}\left(1-\frac{1}{\sqrt{z}}\right) = \frac{\sqrt{z} - \lambda}{2\sqrt{z}}\left(\frac{\sqrt{z}-1}{\sqrt{z-2\sqrt{z}}} - 1\right) \\
& > & 0.
\end{eqnarray*}
So $\sigma$ will enter the region $\mathcal{T}$ as $z$ increases and it cannot escape from this region through the curve $S_{2a}$ as $z$ increases. Suppose $\sigma$ touches the curve $S_{2b}$ at $p = (z,w)$. Let $m_1$ and $m_2$ be the slopes of the tangents of the curves $\sigma$ and $S_{2b}$ at $p$. Then similarly we have
\begin{eqnarray*}
m_1 & = & \frac{1}{2} - \frac{\sqrt{z-2\sqrt{z}}}{2\sqrt{z}} -1, \\
m_2 & = & -1 + \frac{\lambda}{2\sqrt{z}} - \frac{\sqrt{z - 2\sqrt{z}}}{2\sqrt{z}} - \frac{\sqrt{z}-\lambda}{2\sqrt{z - 2\sqrt{z}}}\left(1-\frac{1}{\sqrt{z}}\right),
\end{eqnarray*}
and
\[
m_1 - m_2 = \frac{\sqrt{z} - \lambda}{2\sqrt{z}}\left(\frac{\sqrt{z}-1}{\sqrt{z-2\sqrt{z}}} + 1\right) > 0.
\]
So $\sigma$ will enter the region $\mathcal{T}$ as $z$ increases and it cannot escape from it through the curve $S_{2b}$. It follows that if $\sigma$ meets the boundary of $\mathcal{T}$, then it will enter and stay inside this region as $z \To \infty$.
\end{proof}

We would like to study how the trajectories behave at the infinity. Let $u = \sqrt{z}$, then the equation of $z$ and $w$ is equivalent to the following one in $u$ and $w$:
\begin{equation}\label{eqnuwN6}
\frac{\df w}{\df u} = \frac{2u(\lambda - u - w)}{w}.
\end{equation}
The above equation defines the following dynamical system
\begin{eqnarray*}
\dot{u}(t) & = & w \\
\dot{w}(t) & = & 2u(\lambda - u - w).
\end{eqnarray*}
It has $(\lambda, 0)$ as a critical point which is a stable focus(if $\lambda < 2$) or a stable node(if $\lambda \geq 2$). The system of $u,w$ has one more critical point at $(0,0)$. The linearization shows that it is a saddle point.

The rough picture of the behavior of the system at infinity can be seen from the phase portrait on the Poincar\'{e} sphere by the polar blow-up technique, see for example \cite{ALGM}. Namely, we introduce the new coordinates $X$, $Y$ and $Z$ such that
\[
u = \frac{X}{Z}, \quad w=\frac{Y}{Z}.
\]
Then the phase portrait of the new system in $X, Y, Z$ near the equator $X^2 + Y^2 =1$ of the Poincar\'{e} sphere $\sph^2 = \set{(X,Y,Z)\in \Real^3 : X^2 + Y^2 + Z^2 =1}$ characterizes the behavior of the origin system in $u, w$ at infinity. One can see that $(1,-1, 0)$ is a critical point and it is a saddle-node. The separatrix between the two hyperbolic sectors approaches the focus/node point $(\lambda, 0)$(if $\lambda > \frac{1}{2}$) or the saddle point $(0,0)$(if $\lambda = \frac{1}{2}$) in $uw$-plane. If $\lambda < \frac{1}{2}$, it approaches other critical points $(0,\pm 1,0)$ on the equator of the Poincar\'{e} sphere. When $\lambda \geq \frac{1}{2}$ this separatrix corresponds to the unique trajectory $\gamma$ in $uw$-plane such that along $\gamma$ we have
\[
\lim_{t\To \infty} (u(t),w(t)) = (\lambda, 0) \text{ or }(0,0) \quad \lim_{t\To -\infty} u(t) > 0 \quad \text{and} \quad \lim_{t\To -\infty} \frac{u(t)}{w(t)} = -1.
\]
Note that when $\lambda = \frac{1}{2}$, $\gamma$ is the graph of the function $w = -u$ or $w = -\sqrt{z}$.

In the case when $\lambda > \frac{1}{2}$, if the trajectory $\gamma$ has positive $u$ for every $t$, then there is a corresponding trajectory $\tilde{\gamma}$ in the $zw$-plane with positive $z$ approaching $(\xi, 0)$ as $s\To \infty$ and $\lim_{s\To -\infty} \frac{w}{\sqrt{z}} = -1$. One can see that the integral $I$ is unbounded as $z \To \infty$, i.e., the metric given by the trajectory $\tilde{\gamma}$ is complete.

The above rough picture indicates that we can show the existence of a complete Yamabe shrinker for $\lambda > \frac{1}{2}$ in two steps:
\begin{enumerate}
\item There exists a trajectory $\gamma$ along which, the function $w = w(z)$ is asymptotic to $w = - \sqrt{z}$ for $z$ large(or $s \To - \infty)$.
\item When $s \To \infty$, $\gamma$ meets the $z$-axis with $0 < z \leq \xi$, and then approaches $(\xi, 0)$.
\end{enumerate}

We verify the statements above in the following two lemmas.

\begin{lem}\label{lemN6gammaasym}
There exists a positive number $Z_0$ and a trajectory $\gamma$ of the system (\ref{eqnDSzwshrinker}) in $z, w$ such that
\begin{enumerate}
\item $\gamma$ is in the 4th quadrant when $z> Z_0$, and
\item along $\gamma$ we have
\[
\lim_{z\To \infty} \frac{w}{\sqrt{z}} = -1.
\]
\end{enumerate}
\end{lem}
\begin{proof}
Let $x = \frac{1}{u}$ and $y = \frac{u}{w}$, then the equation (\ref{eqnuwN6}) becomes
\begin{equation}\label{eqnxyN6shrinker}
x^2\frac{\df y}{\df x} = 2\lambda x y^3 - 2 y^3-2y^2 - xy
\end{equation}
We consider a formal power series solution to the above equation
\[
y = x^k(a_0 + a_1 x + a_2 x^2 + \cdots)
\]
such that $a_0 < 0$. Comparing the lowest term on both sides yields two cases. In the first case $k = 1$ and the recurrence relation of $a_i$ is given by
\[
a_0 = -1, \quad a_1 = -2, \quad \lambda = 2
\]
and
\[
(i+2)a_{i} = 2\lambda \sum_{p+q+r = i-2} a_p a_q a_r - 2\sum_{p+q=i} a_p a_q - 2\sum_{p+q+r=i-1}a_pa_q a_r, \quad \text{for } i \geq 3.
\]
So this case can only exist when $\lambda = 2$ and then there is a family of formal power series solutions parameterized by $a_2 \in \Real$.

In the second case we have $k = 0$ and the recurrence relation of $a_i$ is given by
\[
a_0 = -1
\]
and
\[
(i+1)a_i = 2\lambda \sum_{p+q+r= i} a_p a_q a_r - 2 \sum_{p+q = i+1}a_p a_q - 2\sum_{p+q+r=i+1}a_p a_q a_r, \quad \text{for } i \geq 0.
\]
Solving $a_{i+1}$ from the above relation yields
\[
2a_{i+1} = -(i+1)a_i + 2\lambda \sum_{p+q+r = i}a_p a_q a_r - 2\sum_{p+q = i+1}a_p a_q - 2 \sum_{p+q+r = i+1}a_p a_q a_r,
\]
where the indices $p, q$ and $r$ in the last two sums are positive integers. It follows that all the coefficients $a_i$ are determined, i.e., the formal power series in this case is unique. We apply \cite[Theorem 33.1]{Wasow} to the equation (\ref{eqnxyN6shrinker}) and we conclude that there exists an $\eps > 0$ such that for all $x$ with $\abs{x} < \eps$ there is at least one solution $y = \phi(x)$ admits the asymptotic expansion
\[
\phi(x) \sim -1 + \sum_{i=1}^{\infty} a_i x^i, \quad x\To 0.
\]
If we restrict $x$ to be positive and translate back to the variables $z$ and $w$, then it says that for large $z> 0$, there is a solution $w = w(z)$ to the equation (\ref{eqnwzshrinker}) that is asymptotic to the function $w = -\sqrt{z}$. So we can choose $Z_0 > 0$ such that $w=w(z) < 0$ for any $z> Z_0$ and the trajectory $\gamma$ is the one that passes a point in the graph of the solution.
\end{proof}

\begin{rem}
The technique of asymptotic expansion for a differential equation is also used by R. Bryant to show the existence of complete rotationally symmetric steady Ricci soliton metric on $\Real^n$, see \cite{Bryant}.
\end{rem}

\begin{lem}\label{lemN6gammacritical}
Suppose $\gamma$ is the trajectory in the previous lemma. Then
\begin{enumerate}
\item it stays inside the half plane $\mathcal{P}$,
\item it stays between $S_1$ and $S_{2a}$ in the 4th quadrant, and
\item it approaches the critical point $(\xi, 0)$ as $s \To \infty$.
\end{enumerate}
Moreover such trajectory with the above properties is unique.
\end{lem}
\begin{proof}
We already knew that if a trajectory passes a point below the $z$-axis and above $S_1$, then it will meet the $z$-axis and then the value of $w$ stays positive. Similarly since $w''(z) < 0$ in the region $\mathcal{T} \subset \mathcal{P}_2$, any trajectory passes a point in $\mathcal{T}$ will stay in this area as $z$ increases and then $w$ decays at least as fast as the linear function $w = -z$ as $z \To \infty$. Since along $\gamma$ the function $w$ is asymptotically like the function $-\sqrt{z}$, the trajectory $\gamma$ stays between $S_1$ and $S_{2a}$.

If $\lambda \geq 2$, then $S_1$ and $S_{2a}$ meet at $(\xi = \lambda^2, 0)$. Since $\dot{z} = w < 0$, $z$ decreases as $s$ increases. There is no critical point other than $(\xi, 0)$ in the region bounded by $S_1$ and $S_{2a}$, so $\gamma$ will approach this point as $s \To \infty$.

Now we assume that $\frac{1}{2} < \lambda < 2$. We consider another curve $S_3$ in the 4th quadrant:
\[
S_3 = \set{(z,h(z)) : h(z) = -\sqrt{z} \text{ with } z> 0 }.
\]
Suppose $\sigma$ is a trajectory of the system (\ref{eqnDSzwshrinker}) and it touches $S_3$ at $(z,w)$. Then the slope of $\sigma$ at the intersection is given by
\[
m = \frac{\lambda - \sqrt{z}}{-\sqrt{z}} - 1 = -\frac{\lambda}{\sqrt{z}}
\]
and thus we have
\[
m - h'(z) = -\frac{\lambda}{\sqrt{z}} +\frac{1}{2\sqrt{z}} = \left(\frac{1}{2}-\lambda\right)\frac{1}{\sqrt{z}}< 0.
\]
It follows that if $\sigma$ lies below the curve $S_3$ for some $z_0 > 0$, then it stays below $S_3$ for all $z \geq z_0$. Next we consider the relative position between $S_3$ and $S_{2a}$. Since
\[
f_1(z) - h(z)=(\lambda - \sqrt{z})\left(\sqrt{z} - \sqrt{z-2\sqrt{z}}\right) + \sqrt{z} \To \lambda - \frac{1}{2} > 0, \text{ as } z\To \infty,
\]
it follows that $S_{2a}$ lies above $S_3$ for $z>0$ large. In particular if a trajectory lies below $S_3$ for some $z_0 >0$, then it will stay below the curve $S_{2a}$ for all large $z$. So the trajectory $\gamma$ stays in the region $E$ bounded by $S_1$, $S_{2a}$, $S_3$ and the $z$-axis. We consider the convergence of $\gamma$ when $s$ tends to $\infty$. It will meet the $z$-axis at $(z_0, 0)$ with $z_0 \geq 0$. We claim that $z_0 > 0$. If not, then $\gamma$ is a trajectory emanating from the origin. Suppose on $\gamma$ we have $w = w(z)$ for small positive values of $z$. From the differential equation (\ref{eqnwzshrinker}) we have the following series expansion
\[
w(z) = -\sqrt{2\lambda} \sqrt{z} + \text{ higher order terms.}
\]
Since $\lambda > \frac{1}{2}$, for small values of $z$ the trajectory $\gamma$ stays outside the region $E$ which contradicts the fact that $\gamma$ lies inside $E$.  From Proposition \ref{propshrinkercriticalptxi} we conclude that $\gamma$ approaches the stable focus $(\xi, 0)$ within the half plane $\mathcal{P}$ as $s\To \infty$.

To finish the proof, we show that such trajectory is unique. Suppose not, then there are two solutions $w_1(z)$ and $w_2(z)$ which define two trajectories $\gamma_1$ and $\gamma_2$ with the stated properties in this lemma. Since they cannot cross each other in the 4th quadrant, we may assume that $0 < w_2(z) < w_1(z)$ for all $z \geq z_0$ with some positive $z_0$. From the differential equation (\ref{eqnwzshrinker}) we have
\[
w_1'(z) - w_2'(z) = \frac{\lambda - \sqrt{z}}{w_1} - \frac{\lambda - \sqrt{z}}{w_2} = -\frac{\lambda -\sqrt{z}}{w_1 w_2}\left(w_1 - w_2\right).
\]
Since both $\gamma_1$ and $\gamma_2$ stay between $S_1$ and $S_{2a}$, the functions $w_1(z)$ and $w_2(z)$ are asymptotic to the function $-\sqrt{z}$ for large $z$. However from the above differential equation and the comparison of first order differential equations, i.e., the Gr\"{o}nwall's lemma, we see that $\log(w_1 - w_2)$ is asymptotic to the function $2\sqrt{z}$ which shows a contradiction.
\end{proof}

From the previous results we obtain

\begin{thm}\label{thmshrinkerN6}
Let $M^6 = \Real\times_{\vphi} N^5$ and $N$ has constant scalar curvature $\bar{R} > \rho\sqrt{10}$. Then there exists a unique complete Yamabe shrinking soliton metric on $M$ with constant $\rho$ that is not isometric to the Riemannian product.
\end{thm}
\begin{proof}
For any $\lambda > \frac{1}{2}$ with $\rho = \frac{1}{2\sqrt{10}}$, Lemma \ref{lemN6gammaasym} and Lemma \ref{lemN6gammacritical} imply that there exists a trajectory $\gamma$ of the dynamical system (\ref{eqnDSzwshrinker}) such that
\begin{enumerate}
\item $z$ is positive along $\gamma$,
\item $\gamma$ approaches $(\xi, 0)$ with $\xi = \lambda^2$ as $s \To \infty$, and
\item it satisfies the following asymptotic properties:
\[
\lim_{s\To -\infty} z = \infty, \quad \lim_{s\To -\infty} \frac{w}{\sqrt{z}} = -1.
\]
\end{enumerate}
Property (1) ensures that $\vphi$ is positive everywhere, i.e., $\vphi$ defines a Riemannian metric $g = \df r^2 + \vphi^2(r) g_N$ on $M^6 = \Real \times_{\vphi} N^5$. Since the Riemannian distance $r$ is comparable with $s$ near the point $(\xi, 0)$, the metric $g$ is complete when $s\To \infty$. Using the asymptotic properties of $\gamma$ we have
\[
I = \int_{a}^{\infty} \frac{1}{w z^{\frac{1}{4}}}\df z \sim  - \int_{a}^{\infty}z^{-\frac{3}{4}}\df z = -\infty,
\]
for any $a > 0$, i.e., the metric $g$ is complete as $s\To -\infty$. So $\gamma$ defines a complete Yamabe shrinker on $M$. Under the re-scaling of the metric $g$, $\bar{R}/\rho = \lambda/\rho$ is unchanged. So we have the scalar curvature $\bar{R}$ of $N^5$ bigger than $\rho\sqrt{10}$.
\end{proof}

\begin{rem}
From the proof in Lemma \ref{lemN6gammacritical} one sees that the trajectory $\gamma$ is bounded by the curve $w = -\sqrt{z}$ from below. It implies that the scalar curvature $R$ of the metric on $M$ defined by $\gamma$ is positive. To see this, let $\rho = \frac{1}{\sqrt{40}}$ and then we have
\[
R = \rho + \vphi' = \rho + y = \rho + \frac{w}{\sqrt{40}\sqrt{z}} = \frac{1}{\sqrt{40}}\left(1 + \frac{w}{\sqrt{z}}\right) > 0.
\]
\end{rem}

\smallskip

\subsection{The other case with $\dim M \ne 6$} In this case we have
\[
\Phi(z)= \lambda z^{\frac{n-6}{n+2}} - z^{\frac{n-2}{n+2}}.
\]
The two curves that will be useful to study the trajectories of the dynamical system (\ref{eqnDSzwshrinker}) in $z$ and $w$ are
\[
S_1 = \set{(z, w)\in \mathcal{P} : w = \Phi(z)}
\]
and
\[
S_2 = \set{(z,w) \in \mathcal{P} : w^2 + \frac{\Phi(z)}{\Phi'(z)}w - \frac{\Phi^2(z)}{\Phi'(z)} = 0}.
\]
The half plane $\mathcal{P}$ is separated by $S_1$ into two parts:
\[
\mathcal{P}_1 = \set{(z,w) \in \mathcal{P} : w > \Phi(z)}, \quad \mathcal{P}_2 = \set{(z,w) \in \mathcal{P} : w < \Phi(z)}.
\]

When $n < 6$ let $z_{\alpha}$ be the unique positive solution to the equation $4\Phi'(z) + 1 =0$, or equivalently
\[
\frac{n-2}{n-6}z_{\alpha}^{\frac{4}{n+2}} - \frac{n+2}{4(n-6)} z_{\alpha}^{\frac{8}{n+2}} = \lambda.
\]
When $n > 6$ and $(n+2)(n-6)\lambda < (n-2)^2$, the above equation has two distinct positive solutions and let $z_{\alpha}$ be the larger one. If $n > 6$ and $(n+2)(n-6)\lambda \geq (n-2)^2$, then the above equation has no real solution and we let $z_{\alpha} = 0$. It is easy to see that in this case $(\xi, 0)$ is a stable node. In these three cases, we define two functions for $z\geq \max\set{z_{\alpha}, \xi}$:
\begin{eqnarray*}
f_1(z) & = & - \frac{\Phi(z)}{2\Phi'(z)} + \frac{\Phi(z)}{2 \Phi'(z)}\sqrt{1+ 4\Phi'(z)} \\
f_2(z) & = & - \frac{\Phi(z)}{2\Phi'(z)} - \frac{\Phi(z)}{2 \Phi'(z)}\sqrt{1+ 4\Phi'(z)}.
\end{eqnarray*}
We will see that both functions are non-positive and $f_1(z) \geq f_2(z)$ for any $z \geq \max\set{z_\alpha, \xi}$. We call the graphs of $f_{1,2}$ the curves $S_{2a}$ and $S_{2b}$, i.e.,
\begin{eqnarray*}
S_{2a} & = & \set{(z,w)\in \mathcal{P} : w = f_1(z), z \geq \max(z_{\alpha}, \xi)} \\
S_{2b} & = & \set{(z,w)\in \mathcal{P} : w = f_2(z), z \geq \max(z_{\alpha}, \xi)}.
\end{eqnarray*}

\begin{prop}
Let $n < 6$. For these two curves $S_1$ and $S_2$ we have
\begin{enumerate}
\item If $(n+2)\lambda < 16$, i.e., $(\xi, 0)$ is a stable focus, then on $S_2$ we have $z\geq z_\alpha$ and the whole curve $S_2$ is in the 4th quadrant with $S_2 = S_{2a}\cup S_{2b}$ and it is below $S_1$.
\item If $(n+2)\lambda \geq 16$, i.e.,$(\xi, 0)$ is a stable node, then $S_2$ has non-empty components in the 1st quadrant with $z_{\alpha}\leq z < \xi$ and in the 4th quadrant with $z > \xi$. In the 1st quadrant $S_2$ is above $S_1$, and in the 4th quadrant $S_2 = S_{2a}\cup S_{2b}$ and is below $S_1$. Moreover the three curves $S_1$, $S_{2a}$ and $S_{2b}$ have the unique common point $(\xi, 0)$
\end{enumerate}
\end{prop}

\begin{proof}
The defining equation of $S_2$ is given by $w''(z) = 0$ and we have
\[
w''(z) = \frac{\Phi'(z)}{w^2}\left(w^2 + \frac{\Phi(z)}{\Phi'(z)} w - \frac{\Phi^2(z)}{\Phi'(z)}\right) = \frac{\Phi'(z)}{w^2}\left(\left(w + \frac{\Phi(z)}{2\Phi'(z)}\right)^2 - \frac{\Phi^2(z)(1+ 4\Phi'(z))}{4(\Phi'(z))^2}\right).
\]
So on the curve $S_2$ we have $1 + 4\Phi'(z) \geq 0$. We consider the monotonicity of $\Phi'(z)$. Since
\[
\Phi'(z) = \frac{1}{n+2}z^{-\frac{8}{n+2}}\left((n-6)\lambda - (n-2)z^{\frac{4}{n+2}}\right)
\]
and
\[
\Phi''(z) = \frac{1}{(n+2)^2}z^{-\frac{n+10}{n+2}}\left(4(n-2)z^{\frac{4}{n+2}} - 8(n-6)\lambda\right) > 0 \quad \text{for } z > 0,
\]
the function $\Phi'(z)$ is an increasing negative function from $-\frac{1}{4}$ and bounded above by zero as it approaches zero when $z\To \infty$. It follows that $\sqrt{1+ 4\Phi'(z)}$ is an increasing non-negative function which is bounded above by $1$. Let $z_{\alpha} > 0$ be the unique value of $z$ such that $\Phi'(z_{\alpha}) = -\frac{1}{4}$ and then we have
\[
\frac{n-2}{n-6}z^{\frac{4}{n+2}} - \frac{n+2}{4(n-6)}z^{\frac{8}{n+2}} \geq \lambda \quad \text{for } z \geq z_{\alpha},
\]
with equality if $z=z_{\alpha}$. So the defining equations of $S_2$ can be written as
\begin{eqnarray*}
w_1(z) & = & - \frac{\Phi(z)}{2\Phi'(z)} - \frac{\abs{\Phi(z)}}{2 \Phi'(z)}\sqrt{1+ 4\Phi'(z)} \\
w_2(z) & = & - \frac{\Phi(z)}{2\Phi'(z)} + \frac{\abs{\Phi(z)}}{2 \Phi'(z)}\sqrt{1+ 4\Phi'(z)},
\end{eqnarray*}
and the domain of $w_{1,2}$ is $[z_{\alpha}, \infty)$. If $w_1$ or $w_2$ has positive value for some $z \geq z_{\alpha}$, then $\Phi(z) > 0$. It follows that $z_{\alpha} < \xi$ which is equivalent to the following inequality
\[
\frac{n-2}{n-6}\lambda - \frac{n+2}{4(n-6)}\lambda^2 > \lambda,
\]
i.e., $(n+2)\lambda > 16$.

Next we consider the relative position between $S_1$ and $S_2$. Note that when $z\geq \max\set{z_{\alpha},\xi}$ we have $w_{1,2} = f_{1,2}$. If $(n+2)\lambda < 16$, then $S_2$ is in the 4th quadrant and $f_1 \geq f_2$. Since $z_{\alpha} > \xi$, $\Phi(z) < 0$ for $z> z_{\alpha}$ and then we have
\[
\frac{f_1}{\Phi(z)} = -\frac{1}{2\Phi'(z)}\left(1- \sqrt{1+4\Phi'(z)}\right) = \frac{2}{1+ \sqrt{1+ 4\Phi'(z)}} > 1.
\]
So the curve $S_1$ is above $S_2$ in the 4th quadrant. If $(n+2)\lambda = 16$, then a similar argument shows that the two curves $S_1$ and $S_{2}$ have a unique common point $(\xi, 0)$. If $(n+2)\lambda > 16$ then we already know that $S_1$ is above the curve $S_{2}$ in the 4th quadrant. We consider the curves in the 1st quadrant where $z_{\alpha}\leq z < \xi$. Since $\Phi(z) > 0$, we have
\[
\frac{f_{1,2}}{\Phi(z)} = - \frac{1}{2\Phi'(z)}\left(1\pm \sqrt{1+4\Phi'(z)}\right) = \frac{2}{1\mp \sqrt{1+ 4\Phi'(z)}} > 1,
\]
i.e., $S_2$ is above the curve $S_1$. Note that in this case $S_1$ and $S_2$ also have a unique common point $(\xi, 0)$.
\end{proof}

In the case when $n > 6$, if $(n+2)(n-6)\lambda < (n-2)^2$ then $S_2$ has two connected components in the half plane $\mathcal{P}$. The one given by $S_{2a}\cup S_{2b}$ lies entirely in the 4th quadrant. If $(n+2)(n-6)\lambda \geq (n-2)^2$, then $S_2$ has only one component and intersects with $S_1$ at $(\xi, 0)$. Furthermore we have
\begin{prop}
Let $n > 6$. For these two curves $S_1$ and $S_2$ we have
\begin{enumerate}
\item If $(n+2)\lambda < 16$, i.e., $(\xi, 0)$ is a stable focus, then the restriction of $S_2$ to $\set{z \geq \max(z_\alpha, \xi)}$ is $S_{2a}\cup S_{2b}$ which lies entirely in the 4th quadrant and is below $S_1$.
\item If $(n+2)\lambda \geq 16$, i.e.,$(\xi, 0)$ is a stable node, then the restriction of $S_2$ to $\set{z \geq \xi}$ is $S_{2a}\cup S_{2b}$ which lies entirely in the 4th quadrant and is below $S_1$. Moreover the three curves $S_1$, $S_{2a}$ and $S_{2b}$ have the unique common point $(\xi, 0)$.
\end{enumerate}
\end{prop}

\begin{proof}
The argument is similar to the one in the previous proposition. Note that when $z \geq \max\set{z_{\alpha}, \xi}$, we also have $\Phi''(z) > 0$.
\end{proof}

From the formulas of $f_{1,2}$ we see that $f_1 = f_2$ when $z = \max\set{z_{\alpha}, \xi}$ and $f_1$(or $f_2$) is asymptotic to the function $\Phi(z)$(or $-\frac{n+2}{n-2}z$). We define the following trapping area in $\mathcal{P}_2$:
\[
\mathcal{T} = \set{(z,w) \in \mathcal{P} : f_2 \leq w \leq f_1 \text{ and } z \geq \max(z_{\alpha}, \xi)}.
\]

\begin{prop}\label{propshrinkerNgeneralT}
Suppose $\gamma$ is a trajectory of the system (\ref{eqnDSzwshrinker}) with $n \ne 6$. If $\gamma$ meets the boundary of $\mathcal{T}$, then it enters and stays in this area as $z$ increases.
\end{prop}
\begin{proof}
Suppose $\gamma$ meets the boundary of $\mathcal{T}$ at $(z, w)$. Let $m_1$ and $m_2$ be the slopes of the tangents of $\gamma$ and the boundary curve at this point. If $\gamma$ meets $S_{2a}$, then we have
\[
m_1 = \frac{\Phi(z)}{w} - 1 = \frac{\Phi(z)}{f_1(z)} - 1,
\]
and $m_2 = f_1'(z)$. A straightforward computation shows that
\begin{eqnarray*}
m_1-m_2 & = & \frac{\Phi''(z)}{(\Phi'(z))^2\sqrt{1+4\Phi'(z)}}\frac{\Phi(z)}{\sqrt{1+4\Phi'(z)}-1}\left(\sqrt{1+4\Phi'(z)}-1 + \Phi'(z) (\sqrt{1+4\Phi'(z)}-3)\right) \\
& < & 0.
\end{eqnarray*}
So the trajectory $\gamma$ cannot escape the area $\mathcal{T}$ through the curve $S_{2a}$. Similarly if $\gamma$ meets $S_{2b}$ then we have
\begin{eqnarray*}
m_1 - m_2 & = & \frac{\Phi(z)}{f_2(z)} -1 - f_2'(z) \\
& = & - \frac{\Phi''(z)}{(\Phi'(z))^2\sqrt{1+4\Phi'(z)}}\frac{\Phi(z)}{\sqrt{1+4\Phi'(z)}+1}\left(\sqrt{1+4\Phi'(z)}+1 + \Phi'(z) (\sqrt{1+4\Phi'(z)}+3)\right) \\
& > & 0
\end{eqnarray*}
i.e., $\gamma$ cannot leave this area through the curve $S_{2b}$ either. So $\gamma$ will stay inside $\mathcal{T}$ as $z$ increases.
\end{proof}

We have the following existence result of Yamabe shrinkers.
\begin{thm}\label{thmshrinkerNgeneral}
Let $n \ne 6$ and $\rho > 0$. Let $\bar{R} > \rho \frac{n-2}{n+2}\left((n-1)(n+2)\right)^{\frac{4}{n+2}}$ and $(N^{n-1}, \bar{g})$ be Riemannian manifold with constant scalar curvature $\bar{R}$.  Then there exists a unique complete Yamabe shrinking soliton metric on $M^n = \Real \times_{\vphi} N^{n-1}$ with constant $\rho$ that is not isometric to the Riemannian product.
\end{thm}
\begin{proof}
Recall that $\lambda = \bar{R} \left((n-1)(n+2)\right)^{\frac{n-6}{n+2}}$ and $\rho = \left((n-1)(n+2)\right)^{-\frac{n-2}{n+2}}$. The inequality of $\bar{R}$ and $\rho$ in the theorem is equivalent to $\lambda > \frac{n-2}{n+2}$ by re-scaling. As in the case when $n=6$ we show the existence and uniqueness of the metric in two steps.

\textsc{Step 1.} We show that there is a solution $w = w(z)$ for large $z> 0$ with the asymptotic expansion:
\[
w(z)\sim - z^{\frac{n-2}{n+2}} + \text{ lower order terms.}
\]
Then we let $\gamma$ be the trajectory that passes a point in the graph of this solution.

Using the variables $x = z^{-\frac{1}{n+2}}$ and $y = \frac{1}{w}$, the equation (\ref{eqnwzshrinker}) has the following form
\[
x^{2n+1}\frac{\df y}{\df x}= (n+2)\left((\lambda x^4 -1)y^3 - x^{n-2}y^2\right).
\]
For some $\eps>0$ small it has a solution $y = \phi(x)$ with $0< x < \eps$ with the following asymptotic expansion
\[
\phi(x) \sim - x^{n-2} + x^{n-2}\sum_{i=1}^{\infty} a_i x^i, \quad x\To 0.
\]
It follows that the equation (\ref{eqnwzshrinker}) admits a solution $w = w(z)$ such that
\[
\lim_{z\To \infty} \frac{w(z)}{z^{\frac{n-2}{n+2}}} = -1.
\]

\smallskip

\textsc{Step 2.} We show that the trajectory $\gamma$ obtained in the previous step, i.e., it is defined by the solution $w = w(z)$, is unique, stays inside the half plane $\mathcal{P}$ and approaches to the critical point $(\xi, 0)$.

From Proposition \ref{propshrinkerNgeneralT}, $\gamma$ stays between the curve $S_1$ and $S_{2a}$. In the case when $(\xi, 0)$ is a node, i.e., $(n+2)\lambda \geq 16$, since $S_1$ and $S_{2a}$ meet at $(\xi, 0)$, $\gamma$ approaches this point as $z$ decreases. The metric defined by $\gamma$ is complete at $z = \xi$ and the completeness as $z \To \infty$ follows from the asymptotic property of the solution $w = w(z)$.

In the case when $(\xi, 0)$ is a stable focus, we consider the auxiliary curve $S_3$ in the 4th quadrant which is the graph of the function $h(z) = - z^{\frac{n-2}{n+2}}$. One can show the following properties and the proof is similar to the case when $n=6$:
\begin{enumerate}
\item If a trajectory $\sigma$ lies below $S_3$ for some $z_0> 0$, then it will stay below $S_3$ for all $z\geq z_0$.
\item The curve $S_{2a}$ stays below the curve $S_3$ for all large $z> 0$.
\item If a trajectory emanates from the origin, then it will stay below $S_3$.
\end{enumerate}
It follows that the trajectory $\gamma$ stays inside the region $E$ bounded by the curves $S_{1}$, $S_{2a}$, $S_3$ and the $z$-axis, and it meets the $z$-axis with positive value of $z$. From Proposition \ref{propshrinkercriticalptxi} the trajectory $\gamma$ stays inside the half plane $\mathcal{P}$ and approaches the stable focus $(\xi, 0)$. So $\gamma$ defines a complete Yamabe shrinking metric on $M$ which is not isometric to the Riemannian product. The uniqueness of such metric, or equivalently of the trajectory $\gamma$, follows from a similar argument as in the $n=6$ case.
\end{proof}

\begin{rem}
The fact that the scalar curvature $R$ of the metric defined on $M$ by the trajectory $\gamma$ is positive can be seen from the proof in Theorem \ref{thmshrinkerNgeneral}. Let $\rho = \left((n-1)(n+2)\right)^{-\frac{n-2}{n+2}}$ and then we have
\[
R = \rho + y = \rho + w\left((n-1)(n+2) z\right)^{-\frac{n-2}{n+2}} = \rho\left( 1 + w z^{-\frac{n-2}{n+2}}\right) > 0
\]
since the trajectory $\gamma$ is bounded below by the graph of the function $w = - z^{\frac{n-2}{n+2}}$. 
\end{rem}

\medskip



\vfill
\end{document}